\documentclass[twoside,12pt]{article}
\usepackage{amsmath}
\usepackage{amssymb}
\usepackage{amscd}
\usepackage{amsthm}
\theoremstyle{plain}
\newtheorem{Thm}{Theorem}[section]
\newtheorem{Conj}[Thm]{Conjecture}
\newtheorem{Lem}[Thm]{Lemma}
\newtheorem{Prop}[Thm]{Proposition}
\theoremstyle{definition}
\newtheorem{defn}[Thm]{Definition}
\theoremstyle{remark}
\newtheorem{Rem}[Thm]{Remark}
\def\Strata{\mathop{\mathbf{Strata}}}
\def\MS{\mathop{\mathbf{MS}}}

\title{\textbf{A note on the ampleness of numerically positive log canonical and anti-log canonical divisors}}
\author{\thanks{2000 \textit{Mathematics Subject Classification}. 14E30} \textbf{Shigetaka Fukuda}}
\date{\empty}

\begin{document}
\maketitle \thispagestyle{empty}
\pagestyle{myheadings}
\markboth{Shigetaka Fukuda}{A note on the ampleness}
\begin{abstract}
In this short note, we consider the conjecture that the log canonical divisor (resp.\ the anti-log canonical divisor) $K_X + \Delta$ (resp.\ $-(K_X + \Delta)$) on a pair $(X, \Delta)$ consisting of a complex projective manifold $X$ and a reduced simply normal crossing divisor $\Delta$ on $X$ is ample if it is numerically positive.
More precisely, we prove the conjecture for $K_X + \Delta$ with $\Delta \neq 0$ in dimension $4$ and for $-(K_X + \Delta)$ with $\Delta \neq 0$ in dimension $3$ or $4$.
\end{abstract}

Every variety is defined over the field of complex numbers throughout the paper.
Let $X$ be an $n$-dimensional nonsingular projective algebraic variety and $\Delta = \sum_{i \in I} \Delta_i $ a reduced simply normal crossing divisor on $X$ (where $\Delta_i $ is a prime divisor).
We denote the canonical divisor of $X$ by $K_X$.
Thus $K_X + \Delta$ denotes the log canonical divisor on the pair $(X, \Delta)$.

By the symbol $\kappa (X, L)$, we mean the Iitaka dimension of a $\mathbf{Q}$-Cartier $\mathbf{Q}$-divisor $L$.

\begin{defn} A $\mathbf{Q}$-Cartier $\mathbf{Q}$-divisor $L$ on $X$ is {\it numerically positive} ({\it nup} (\cite{LaRo}), for short) if $(L, C) > 0$ for every curve $C$ on $X$.
\end{defn}

\begin{Rem}
In the case where $\dim X = 1$, the nupness means the ampleness.
\end{Rem}

In this paper we deal with the following four conjectures, which are well known to the specialists of higher dimensional algebraic varieties.

\begin{Conj}\label{C1}
If $K_X$ is nup, then it is ample.
\end{Conj}

\begin{Conj}\label{C2}
If $K_X + \Delta$ is nup, then it is ample.
\end{Conj}

\begin{Conj}\label{C3}
If $-K_X$ is nup, then it is ample.
\end{Conj}

\begin{Conj}\label{C4}
If $-(K_X + \Delta)$ is nup, then it is ample.
\end{Conj}

Conjectures \ref{C1} and \ref{C2} are theorems in dimension $n \leq 3$, by virtue of the abundance and the log abundance theorems (due to Kawamata \cite{Ka87} and \cite{Ka92}, Miyaoka \cite{Mi} and Keel-Matsuki-McKernan \cite{KeMaMc}).
Conjecture \ref{C3} was proved by Hidetoshi Maeda \cite{MaHd} in dimension $n=2$ and by Serrano \cite{Se} in dimension $n=3$.

Hironobu Maeda \cite{MaHr} proved Conjecture \ref{C4} in the case where $\Delta \neq 0$ and $n=2$, as follows:
Assume that n=2, the anti-log canonical divisor $-(K_X + \Delta)$ is nup and $\Delta \ne 0$.
First we shall show that $(-(K_X + \Delta))^2 > 0$.
Let us derive a contradiction, assuming that $(-(K_X + \Delta))^2 = 0$.
From the nupness of $-(K_X + \Delta)$, we have $-(K_X + \Delta) \Delta > 0$.
Thus $-(K_X + \Delta) K_X < 0$.
Then $\kappa (X, -(K_X + \Delta)) = 1$ by virtue of Sakai \cite{Sa}, Theorem 2.
Hence the nupness of $-(K_X + \Delta)$ implies that $(-(K_X + \Delta))(-(K_X + \Delta)) > 0$, because a high multiple of $-(K_X + \Delta)$ becomes linearly equivalent to some nonzero effective divisor.
This is a contradiction!
Consequently we have $(-(K_X + \Delta))^2 > 0$.
Next we apply the Nakai criterion to the divisor $-(K_X + \Delta)$ and obtain that it is ample.

By using Wilson's technique \cite{Wi}, Hironobu Maeda \cite{MaHr} proved Conjecture \ref{C4} also in dimension $n=3$ under the extra condition $\kappa (X, -(K_X + \Delta)) \geq 1$.
(This result was reviewed by Matsuki \cite{Ma87}.)

Here we remark that Serrano \cite{Se} has implicitly proved Conjecture \ref{C4} in dimension $n=3$ under the weaker condition that $\kappa (X, -(K_X + \Delta)) \geq 0$, as follows:
Assume that $n=3$, that the anti-log canonical divisor $-(K_X + \Delta)$ is nup and that $\kappa (X, -(K_X + \Delta)) \geq 0$.
Then $\kappa (X, (-1)K_X +1(-(K_X + \Delta))) = \kappa (X, -2(K_X + \Delta) + \Delta) \geq \kappa (X, -2(K_X + \Delta)) = \kappa (X, -(K_X + \Delta)) \geq 0$.
Thus Serrano \cite{Se}, Proposition 3.1 implies that $-(K_X + \Delta) + \epsilon K_X$ is ample for a sufficiently small positive rational number $\epsilon$.
Therefore $-(K_X + \Delta) = (1/(1- \epsilon))((-(K_X + \Delta) + \epsilon K_X )+ \epsilon \Delta)$ is big.
This satisfies the extra condition stated in the preceding paragraph.

Now we state our main theorem

\begin{Thm}\label{MT}
$(1)$ Conjecture \ref{C2} is true in the case where $\Delta \neq 0$ and $n=4$.
$(2)$ Conjecture \ref{C4} is true in the case where $\Delta \neq 0$ and $n=3,4$.
\end{Thm}

{\bf Acknowledgment:}
I thank the referee for careful reading the manuscript and for valuable advice concerning the presentation.

\section{Proof of Theorem \ref{MT}}

We define $\Strata (\Delta) := \{\Gamma \mid \Gamma \text{ is an } \text{irreducible } \text{component } \text{of } \bigcap_{j \in J} \Delta_j \ne \emptyset \text{,} \text{ for } \text{some } \text{nonempty } \text{subset } J \text{ of } I \}$ and $\MS (\Delta) := \{\Gamma \in \Strata (\Delta) \mid \text{ If } \Gamma ' \in \Strata (\Delta) \text{ and } \Gamma ' \subseteq \Gamma, \text{ then } \Gamma ' = \Gamma \}$.
We remark that $(K_X + \Delta) \mid_{\Gamma} = K_{\Gamma}$ for every $\Gamma \in \MS (\Delta)$.

Let $L$ be a $\mathbf{Q}$-Cartier $\mathbf{Q}$-divisor on $X$.

$L$ is said to be {\it nef and log big} on $(X, \Delta)$, if $L$ is nef, $L^n > 0$ and $(L \mid_{\Gamma})^{\dim \Gamma} > 0$ for any $\Gamma \in \Strata (\Delta)$.

\begin{Rem}
Assume that $L$ is nef.

If $bL-(K_X + \Delta)$ is nef for some $b \geq 0$, then so is $aL-(K_X + \Delta)$ for $a \gg 0$.

If $bL-(K_X + \Delta)$ is nup for some $b \geq 0$, then so is $aL-(K_X + \Delta)$ for $a \gg 0$.

If $bL-(K_X + \Delta)$ is nef and big for some $b \geq 0$, then so is $aL-(K_X + \Delta)$ for $a \gg 0$.

If $bL-(K_X + \Delta)$ is nef and log big on $(X, \Delta)$ for some $b \geq 0$, then so is $aL-(K_X + \Delta)$ for $a \gg 0$.
\end{Rem}

We cite two lemmas:

\begin{Lem}[An uniruledness theorem of Miyaoka-Mori type, Matsuki \cite{Ma03}]\label{L1}
Let $D_1, D_2, \ldots , D_n$ be a sequence of nef Cartier divisors.
Suppose $D_1 \cdot D_2 \cdots D_n = 0$ and $-K_X \cdot D_1 \cdot D_2 \cdots D_{n-1} >0$.
Then $X$ is covered by a family of rational curves $C$ such that $D_n \cdot C =0$.
\end{Lem}

\begin{Lem}[Base point free theorem of Reid type, Fukuda \cite{Fu}]\label{L2}
If $L$ is nef and $bL-(K_X + \Delta)$ is nef and log big on $(X, \Delta)$ for some $b \geq 0$, then $L$ is semi-ample.
\end{Lem}

\begin{Prop}\label{P}
Assume that $L$ is nef and $bL-(K_X + \Delta)$ is nup for some $b \geq 0$ and that $\Delta \ne 0$.
If $((bL-(K_X + \Delta)) \mid_{\Gamma})^{\dim \Gamma} > 0$ for any $\Gamma \in \MS (\Delta)$, then $aL-(K_X + \Delta)$ is nef and log big on $(X, \Delta)$ for $a \gg 0$.
\end{Prop}

\begin{proof}
We prove this proposition by induction on $n$.
If $n=1$, the statement is trivial.
Thus we may assume that $n \ge 2$.

We note that $(aL-(K_X + \Delta)) \mid_{\Delta_i} = aL \mid_{\Delta_i} -(K_{\Delta_i} + (\Delta - \Delta_i) \mid_{\Delta_i})$.

First we shall show that $((aL-(K_X + \Delta)) \mid_{\Gamma})^{\dim \Gamma} > 0$ for any $\Gamma \in \Strata (\Delta)$.
If $(\Delta - \Delta_i) \mid_{\Delta_i} \ne 0$, then the induction hypothesis implies that $((aL-(K_X + \Delta)) \mid_{\Gamma})^{\dim \Gamma} > 0$ for any $\Gamma \subseteq \Delta_i $.
Thus we may assume that $(\Delta - \Delta_i) \mid_{\Delta_i} = 0$.
Then $\Delta_i \in \MS (\Delta)$.
Therefore $((aL-(K_X + \Delta)) \mid_{\Delta_i})^{\dim \Delta_i} > 0$.

Next we shall show that $(aL-(K_X + \Delta))^n > 0$.
Assuming that $(aL-(K_X + \Delta))^n = 0$ for any $a \gg 0$, we will derive the contradiction.
Then we have $L^i (K_X + \Delta)^{n-i} = 0$ for $i = 0,1,2, \ldots , n $, by regarding $(aL-(K_X + \Delta))^n$ as a polynomial in the variable $a$.
Thus $-K_X \cdot (aL-(K_X + \Delta))^{n-1} = (-aL + \Delta) \cdot (aL-(K_X + \Delta))^{n-1} = \Delta \cdot (aL-(K_X + \Delta))^{n-1} \ge (aL-(K_X + \Delta))^{n-1} \Delta_i = ((aL-(K_X + \Delta)) \mid_{\Delta_i})^{\dim \Delta_i} > 0$.
Consequently Lemma \ref{L1} derives the contradiction.
\end{proof}

\begin{Thm}\label{STH}
Assume that $L$ is nef and $bL-(K_X + \Delta)$ is nup for some $b \geq 0$ and that $\Delta \ne 0$.
If $((bL-(K_X + \Delta)) \mid_{\Gamma})^{\dim \Gamma} > 0$ for any $\Gamma \in \MS (\Delta)$, then $L$ is semi-ample.
\end{Thm}

\begin{proof}
The assertion follows from Lemma \ref{L2}, because $aL-(K_X + \Delta)$ is nef and log big on $(X, \Delta)$ for $a \gg 0$ by Proposition \ref{P}.
\end{proof}

\begin{Prop}\label{P2}
$(1)$ Conjecture \ref{C2} is true in the case $\Delta \ne 0$, if Conjecture \ref{C1} is true in dimension $\le n-1$.
$(2)$ Conjecture \ref{C4} is true in the case $\Delta \ne 0$, if Conjecture \ref{C3} is true in dimension $\le n-1$.
\end{Prop}

\begin{proof}
(1) Put $L = K_X + \Delta$ in the statement of Theorem \ref{STH}.
(2) Put $L = -(K_X + \Delta)$ in the statement of Theorem \ref{STH}.
\end{proof}

\begin{proof}[Proof of Theorem \ref{MT}]
(1): Conjecture \ref{C1} is true in the case $n \le 3$ (Miyaoka \cite{Mi}, Kawamata \cite{Ka92}).
Thus Proposition \ref{P2} implies the assertion.

(2): Conjecture \ref{C3} is true in the case $n \le 3$ (Hidetoshi Maeda \cite{MaHd}, Serrano \cite{Se}).
Thus Proposition \ref{P2} implies the assertion.
\end{proof}

\bigskip
Faculty of Education, Gifu Shotoku Gakuen University

Yanaizu-cho, Gifu 501-6194, Japan

fukuda@ha.shotoku.ac.jp

\end{document}